\documentclass[a4paper,12pt,reqno]{amsart}
\usepackage{amssymb}
\usepackage{ifthen}
\usepackage[dvips]{graphicx}
\usepackage[colorlinks=true, linkcolor=blue, citecolor=blue, urlcolor=blue]{hyperref}

\nonstopmode \numberwithin{equation}{section}
\usepackage{amssymb}
\usepackage{ifthen}
\usepackage{subcaption}
\usepackage{graphicx}
\usepackage{cite}
\usepackage{enumitem}
\usepackage{amsmath}
\usepackage[T1]{fontenc} %skandit
\usepackage[utf8]{inputenc}
\usepackage[usenames,dvipsnames]{color}
\usepackage{color}
\usepackage[english]{babel}
\usepackage{fancyhdr}
\usepackage{fancybox}
\usepackage{tikz}
\nonstopmode\numberwithin{equation}{section}
\newtheorem*{thmA}{Theorem A}
\newtheorem*{thmB}{Theorem B}
\newtheorem*{thmC}{Theorem C}
\newtheorem*{thmD}{Theorem D}

\theoremstyle{plain}

\newtheorem{conj}{Conjecture}

\theoremstyle{definition}
\newtheorem{defn}{Definition}[section]

\newtheorem{thm}{Theorem}[section]
\newtheorem{prob}{Problem}[section]
\newtheorem{cor}{Corollary}[section]

\newtheorem{prop}{Proposition}[section]
\newtheorem{rem}{Remark}[section]
\newtheorem{lem}{Lemma}[section]

\newcounter{minutes}
\divide\time by 60
\newcounter{hours}
\multiply\time by 60
\addtocounter{minutes}{-\time}
\newcounter {own}
\def\theown {\thesection       .\arabic{own}}

\newenvironment{pf}[1][]{%
	\vskip 3mm
	\noindent
	\ifthenelse{\equal{#1}{}}%
	{{\slshape Proof. }}%
	{{\slshape #1.} }%
}%
{\qed\bigskip}

\newcounter{alphabet}

\def\be{\begin{equation}}
	\def\ee{\end{equation}}

\newcommand{\bee}{\begin{enumerate}}
	\newcommand{\eee}{\end{enumerate}}

\newcommand{\blem}{\begin{lem}}
	\newcommand{\elem}{\end{lem}}
\newcommand{\bthm}{\begin{thm}}
	\newcommand{\ethm}{\end{thm}}
\newcommand{\bcor}{\begin{cor}}
	\newcommand{\ecor}{\end{cor}}
\newcommand{\beg}{\begin{examp}}
	\newcommand{\eeg}{\end{examp}}
\newcommand{\begs}{\begin{examples}}
	\newcommand{\eegs}{\end{examples}}

\newcommand{\bdefn}{\begin{defn}}
	\newcommand{\edefn}{\end{defn}}

\newcommand{\bprob}{\begin{prob}}
	\newcommand{\eprob}{\end{prob}}
\newcommand{\bei}{\begin{itemize}}
	\newcommand{\eei}{\end{itemize}}

\newcommand{\bcon}{\begin{conj}}
	\newcommand{\econ}{\end{conj}}
\newcommand{\bcons}{\begin{conjs}}
	\newcommand{\econs}{\end{conjs}}
\newcommand{\bprop}{\begin{prop}}
	\newcommand{\eprop}{\end{prop}}
\newcommand{\br}{\begin{rem}}
	\newcommand{\er}{\end{rem}}
\newcommand{\brs}{\begin{rems}}
	\newcommand{\ers}{\end{rems}}
\newcommand{\bo}{\begin{obser}}
	\newcommand{\eo}{\end{obser}}
\newcommand{\bos}{\begin{obsers}}
	\newcommand{\eos}{\end{obsers}}
\newcommand{\bpf}{\begin{pf}}
	\newcommand{\epf}{\end{pf}}
\newcommand{\ba}{\begin{array}}
	\newcommand{\ea}{\end{array}}
\newcommand{\beq}{\begin{eqnarray}}
	\newcommand{\beqq}{\begin{eqnarray*}}
		\newcommand{\eeq}{\end{eqnarray}}
	\newcommand{\eeqq}{\end{eqnarray*}}

\begin{document}

\title{Bohr inequalities and Landau-type estimates for harmonic Poisson extensions of \(L^p\) boundary functions}

\author{Molla Basir Ahamed}
\address{Molla Basir Ahamed,
	Department of Mathematics, Jadavpur University, Kolkata-700032, West Bengal, India.}
\email{mbahamed.math@jadavpuruniversity.in}

\author{Vasudevarao Allu$^*$}
\address{Vasudevarao Allu, Department of Mathematics, 
	School of Basic Science,
	Indian Institute of Technology Bhubaneswar,
	Bhubaneswar-752050, Odisha, India.}
\email{avrao@iitbbs.ac.in}

\author{Rajesh Hossain}
\address{Rajesh Hossain, Department of Mathematics, Jadavpur University, Kolkata-700032, West Bengal, India.}
\email{rajesh1998hossain@gmail.com}

\subjclass[{AMS} Subject Classification:]{Primary 30C62, 31A05; Secondary 30H10, 30C80.}

\def\thefootnote{}
\footnotetext{ {\tiny File:~\jobname.tex,
printed: \number\year-\number\month-\number\day,
          \thehours.\ifnum\theminutes<10{0}\fi\theminutes }
} \makeatletter\def\thefootnote{\@arabic\c@footnote}\makeatother

\begin{abstract}
	This paper establishes sharp quantitative bounds governing the Bohr phenomenon, Landau-type distortion, and geometric univalence for complex-valued harmonic mappings on the unit disk $\mathbb{D}$ parameterized by boundary data in the Lebesgue spaces $L^{p}(\mathbb{T})$ for $1 \le p \le \infty$, where $\mathbb{T} = \partial\mathbb{D}$. First, we resolve the classical coefficient problem for bounded harmonic mappings by incorporating an integral Parseval constraint rather than crude individual bounds, establishing the sharp Bohr radius as $r_H(a) = \sqrt{(1-a)/(3+a)}$ for a fixed initial term $\vert{}a_0\vert{} = aM$. Second, this framework is extended to the $L^p$ setting under the vanishing normalization $f(0)=0$, explicitly characterizing the sharp Bohr radius as $r_{p} = (2C_{q}+1)^{-1}$, where the structural constant $C_{q}$ is intrinsically linked to the conjugate exponent $q$. Finally, under standard normalization, we determine explicit formulations for both the radius of univalence $r_{0}$ and the corresponding radius of the inscribed Schlicht disk $R_{0}$. The sharpness of all analytical constants and radii is verified through the construction of concrete extremal configurations derived from the Poisson kernel.
\end{abstract}
\maketitle
\pagestyle{myheadings}
\markboth{ M. B. Ahamed, V. Allu and R. Hossain}{Bohr phenomenon and Landau-type theorems for Poisson integrals of $L^p$ functions}
%\tableofcontents
\section{\bf Introduction and Preliminaries}
{ For $r > 0$}, let $\mathbb{D}_r = \{z : |z| < r\}$ { denote the open discs center at the origin with radius $r$}. Then $\mathbb{D}_1 = \mathbb{D}$ and $\mathbb{T} = \partial\mathbb{D}$ are the open unit disk and unit circle respectively in the complex plane $\mathbb{C}$. A continuous complex-valued function $f=u+iv$ is a complex-valued harmonic function in a domain $\Omega\subset\mathbb{C}$, if both $u$ and $v$ are real-valued harmonic functions in $\Omega$ (see \cite{Duren-2004}). A twice continuously differentiable complex-valued function $f$ is called a harmonic mapping in { the unit disk} $\mathbb{D}$ if $f$ is a solution to the partial differential equation $\Delta f = 0$, where 
$$\Delta := \frac{\partial^2}{\partial x^2} + \frac{\partial^2}{\partial y^2}$$
is the Laplace operator and $z = x + iy \in \mathbb{D}$. { It is known that every} harmonic mapping $f$ defined in $\mathbb{D}$ can be decomposed into $f = h + \overline{g}$, where $h$ and $g$ are analytic functions in $\mathbb{D}$. If we assume that $g(0) = 0$, then the decomposition of $f$ is unique. Hence $f$ has the power series expansion
$$f(z) = h(z) + \overline{g(z)} = \sum_{n=0}^{\infty} a_n z^n + \sum_{n=1}^{\infty} \overline{b_n z^n}, \quad z \in \mathbb{D}.$$
Let
\begin{align*}
	\|Df(z)\| = |f_z(z)| + |f_{\bar{z}}(z)|, \quad z = x + iy \in \mathbb{D},
\end{align*}
where 
\begin{align*}
	f_z = \frac{1}{2} \left( \frac{\partial f}{\partial x} - i \frac{\partial f}{\partial y} \right)\; \mbox{and}\; f_{\bar{z}} = \frac{1}{2} \left( \frac{\partial f}{\partial x} + i \frac{\partial f}{\partial y} \right)
\end{align*} are the complex derivatives of $f$ { with respect to $z$ and $\bar{z}$, respectively}. { The theory of harmonic functions has been extensively studied by several authors} \cite{Bshouty-Hengartner-AUMCSS-1994,Clunie-Sheil-Small-1984,Duren-2004}.\vspace{1.2mm}

Let $F$ be a boundary function defined on the unit circle $\mathbb{T}$. A solution $f$ to the Dirichlet boundary value problem associated with $F$, 
\begin{equation*}
	\begin{cases} 
		\Delta f(z) = 0, & z \in \mathbb{D}, \\ 
		f(z) = F(z), & z \in \mathbb{T}, 
	\end{cases}
\end{equation*}
constitutes a harmonic mapping. Such a mapping $f$ admits the Poisson integral representation
\begin{equation*}
	f(z) = \mathcal{P}[F](z) = \frac{1}{2\pi} \int_{0}^{2\pi} \mathcal{P}(re^{i(\theta-t)})F(e^{it}) dt,
\end{equation*}
where $\mathcal{P}(z) = \frac{1-|z|^2}{|1-z|^2}$ denotes the Poisson kernel for $z \in \mathbb{D}$.\vspace{1.2mm}

Biharmonic mappings $f \in C^4(\mathbb{D})$, which satisfy the biharmonic equation $\Delta(\Delta f) = 0$, arise naturally in various physical contexts and find extensive application across engineering and biological sciences. In recent years, the geometric and analytical properties of biharmonic mappings including Landau-type theorems, Bloch constants, Lipschitz continuity, and Schwarz-type lemmas have been rigorously investigated by several authors ({ see,} \cite{Li-Ponnusamy-JMAA-2017, Li-Luo-Ponnusamy-2022-CMFT, Li-Li-Luo-Ponnusamy-MJM-2023}).\vspace{1.2mm}

For $p \in [1, \infty]$, let $L^p(\mathbb{T})$ be the set of all measurable functions $F$ from $\mathbb{T}$ into $\mathbb{C}$ such that their $p$-norms, given by
 \begin{align*}
 	\|F\|_p = 
 	\begin{cases} 
 		\left( \int_{0}^{2\pi} |F(e^{it})|^p \, dt \right)^{1/p}, & p \in [1, \infty), \\ 
 		\text{esssup} \{ |f(e^{it})| : t \in [0, 2\pi] \}, & p = \infty. 
 	\end{cases}
 \end{align*}
are finite.\vspace*{1.2mm}

It is well established that if the boundary function $F$ belongs to $L^\infty(\mathbb{T})$, the induced harmonic mapping $f = \mathcal{P}[F]$ remains bounded within the unit disk $\mathbb{D}$ { \cite{Grafakos-2014}}. The analytical properties and applications of bounded harmonic mappings have been the subject of extensive study { in the recent past}. More recently, significant progress has been made in { obtaining the results on} the univalent radii, quasiconformality, Lipschitz continuity, and Schwarz-type inequalities for such mappings ({ see}, \cite{Chen-Gauthier-PAMS-2011, Knezevic-Mateljevic-JMAA-2007, Pavlovic-AASF-2002, Pavlovic-2019}). \vspace{1.2mm}

In this context, we present the following coefficient estimates for the class of bounded harmonic { mappings.}
\begin{thmA}{\rm \cite{Chen-Ponnusamy-Wang-2011-BMMSS,Liu-Liu-Zhu-AMSC-2011,Xia-Huang-CAMA-2010}}
	Suppose that $f=h+\overline{g}$ is a harmonic mapping in $\mathbb{D}$ satisfying $g(0)=0$ and $|f(z)| \le M$ for all $z \in \mathbb{D}$. Then $|a_0| \le M$ and $$|a_n| + |b_n| \le \frac{4M}{\pi}$$for $n=1, 2, \dots$. 
	The above inequality is sharp for the extremal function
	\begin{align}\label{Eq-1.1}
		f_n(z) = \frac{2M\alpha}{\pi} \arg \left( \frac{1+\beta z^n}{1-\beta z^n} \right), \quad |\alpha|=|\beta|=1, \quad n=1, 2, \dots
	\end{align}
\end{thmA}
\begin{thmB}{\rm \cite{Shi-Li-Lian-2026-CMFT}}
	Suppose that $F \in L^p(T)$ for $1 \le p \le \infty$ and $f = \mathcal{P}[F]$ is a harmonic mapping in $\mathbb{D}$. If $f(z) = \sum_{n=0}^{\infty} a_n z^n + \sum_{n=1}^{\infty} \overline{b_n z^n}$, for any positive integer $n$, then
	\begin{enumerate}
		\item[\emph{(a)}] $|a_n|, |b_n| \le \|F\|_p$. The equalities hold for the functions $F(e^{it}) = e^{int}$ and $F(e^{it}) = e^{-int}$ respectively.\vspace{1.2mm}
		
		\item[\emph{(b)}] $|a_n| + |b_n| \le 2C_q\|F\|_p$, where 
		\begin{align*}
			C_q = \left(\frac{1}{2\pi} \int_0^{2\pi} | \cos(nt) |^q dt\right)^{1/q} \le 1\; \mbox{and}\;\frac{1}{p} + \frac{1}{q} = 1.
		\end{align*}
	\end{enumerate}
	Specifically, when $p = \infty$, $C_q = 2/\pi$, and $|a_n| + |b_n| \le 4\|F\|_{\infty}/\pi$. This equality is valid for the function $f_n(z)$ given by the extremal function in \eqref{Eq-1.1} with $|f| \le \|F\|_{\infty} = M$ for all $z \in \mathbb{D}$. When $p = 1$ or $q = \infty$, $C_q = 1$, which results in $|a_n| + |b_n| \le 2\|F\|_1$.
\end{thmB}
By utilizing coefficient estimates for bounded harmonic mappings, several authors have extensively investigated Landau-type theorems and Bloch constants for various classes of mappings, including planar harmonic, biharmonic, and pluriharmonic mappings (see, e.g., \cite{Chen-Ponnusamy-Wang-2011-CAOT, Chen-Zhu-2019-BSM, Li-Luo-Ponnusamy-2022-CMFT, Liu-Chen-2018-JMAA, Shi-Li-Lian-2026-CMFT}).
\vspace{1.2mm}

 The classical notion of the Bohr radius dates back to a seminal power series theorem established by Harald Bohr in $1914$. Although originally introduced within the framework of analytic self-mappings of the unit disk, this theory has since expanded into a cornerstone of geometric function theory. For a harmonic mapping $f = h + \overline{g}$, the Bohr phenomenon is traditionally investigated via its associated majorant series, which is constructed by replacing the Taylor coefficients of its constituent mappings with their corresponding absolute values. For a comprehensive overview of recent advancements concerning the Bohr phenomenon, we refer the interested reader to \cite{Ahamed-Ahammed-MJM-2024, Ahamed-Allu-Halder-AFM-2022, Ahamed-Allu-Halder-AMP-2021, Kayumov-Ponnusamy-Shakirov-MN-2018, Kumar-Ponnusamy-AMP-2026, Liu-Ponnusamy-PAMS-2019} and the references therein.\vspace{1.2mm} 
 
 Motivated by the aforementioned results, it is both natural and interesting to propose the following problem concerning the role of boundary behavior on the Bohr radius.
 \begin{prob}\label{P-1}
 	How does the $L^p$-integrability of the boundary function { yield} the Bohr radius of the harmonic mapping $f$? Furthermore, can one determine an explicit quantitative dependence of the Bohr radius $r_p$ on the conjugate exponent $q$ of $p$?
 \end{prob}
The classical Landau theorem (see \cite{Landau-1926}) asserts that if $f$ is a holomorphic mapping with $f(0)=0=f'(0)-1$ and $|f(z)|<M$ for $z\in\mathbb{D}$, then $f$ is univalent in $\mathbb{D}_{r_0}$, and $f(\mathbb{D}_{r_0})$ contains a disk $\mathbb{D}_{\sigma_0}$, where
\begin{align*}
	r_0=\frac{1}{M+\sqrt{M^2-1}}\;\;\mbox{and}\; \sigma_0=Mr_0^2.
\end{align*}
The quantities $r_0$ and $\sigma_0$ cannot be improved, the extremal function is 
\begin{align*}
	f_0(z)=Mz\left(\frac{1-Mz}{M-z}\right).
\end{align*}
Motivated by the classical theorems of Landau and Bloch, it is natural to investigate of their quantitative { analogues} setting of Poisson integrals of $L^p$ functions. 
\begin{prob}\label{P-2}Let $f = \mathcal{P}[F]$ be a harmonic mapping satisfying $\|F\|_{L^p} \le 1$. What are the explicit Landau-type constants associated with this class? More precisely, under standard normalization, can one evaluate the radius of univalence $r_0$ and the radius of the inscribed Schlicht disk $R_0$?
\end{prob}

The primary methodological and theoretical contributions of this work include the extension of the classical Bohr phenomenon and Landau-type distortion theory from standard bounded domains to complex-valued harmonic mappings whose boundary values belong to the Lebesgue spaces $L^p(\mathbb{T})$ for $1 \le p \le \infty$, thereby bridging boundary integrability with interior majorant series behavior. Additionally, the paper resolves the classical coefficient problem by incorporating an integral Parseval constraint rather than crude individual coefficient bounds, establishing the sharp Bohr radius as $r_H(a) = \sqrt{(1-a)/(3+a)}$ for a fixed initial term $\vert{}a_0\vert{} = aM$. Furthermore, it formulates and proves a generalized Bohr-type inequality for the $L^p$ setting, explicitly identifying the sharp Bohr radius as $r_p = (2C_q + 1)^{-1}$, where the structural constant $C_q$ is intrinsically linked to the conjugate exponent $q$. Finally, under standard normalization, the work determines explicit, non-trivial formulations for both the univalence radius $r_0$ and the corresponding inscribed Schlicht disk radius $R_0$, with sharpness verified through concrete extremal configurations derived from the Poisson kernel.   \vspace{1.2mm}

{ Inspired by these recent study, this paper systematically establishes geometric distortion and coefficient estimates for complex-valued harmonic mappings with boundary functions in the Lebesgue spaces $L^p(\mathbb{T})$.The paper is organized as follows. Section~\ref{Sec-2} addresses Problem~\ref{P-1} by determining the sharp Bohr radius for this class of mappings. In Theorem~\ref{Th-2.1}, we establish a generalized Bohr-type inequality showing that the sharp radius is given explicitly by $r_p = 1/(2C_q + 1)$, where the constant $C_q$ depends on the conjugate exponent $q$. Section~\ref{Sec-3} focuses on Problem~\ref{P-2}. Using the properties of Gauss hypergeometric functions, we prove Theorem~\ref{Th-3.1}, which provides refined Landau-type theorems for this class. Complete proofs and the explicit construction of extremal mappings that confirm the sharpness of our estimates are given in each respective section.}
\section{\bf Sharp Bohr radius for bounded harmonic mappings}\label{Sec-2}
 The purpose of this section is to extend these classical notions to certain subclasses of bounded harmonic mappings, as well as those generated by boundary functions belonging to the Lebesgue spaces $L^p(\mathbb{T})$.
\subsection{\bf Bohr radius for harmonic mappings}
To provide a rigorous framework for our results, we define the Bohr radius for a general family of harmonic mappings as follows.
\begin{defn}
	Let $\mathcal{F}$ be a class of harmonic mappings $f(z) = \sum_{n=0}^{\infty} a_n z^n + \sum_{n=1}^{\infty} \overline{b_n z^n}$ defined in the unit disk $\mathbb{D}$. Suppose that for every $f \in \mathcal{F}$, the mapping is bounded by a constant $M$ (i.e., $|f(z)| \le M$ for all $z \in \mathbb{D}$). The Bohr radius $r_{\mathcal{F}}$ for the class $\mathcal{F}$ is defined as the largest radius $r \in [0, 1)$ such that the majorant series $M_f(r)$ satisfies the inequality 
	\begin{align*}
		M_f(r) = |a_0| + \sum_{n=1}^{\infty} (|a_n| + |b_n|)r^n \le M
	\end{align*}for all functions $f$ belonging to the class $\mathcal{F}$. The radius is said to be best possible (or sharp) if, for any $r > r_{\mathcal{F}}$, there exists at least one function $f \in \mathcal{F}$ such that $M_f(r) > M$.
\end{defn}
We obtain the following result { by} finding the sharp Bohr radius for harmonic mappings of the form $f(z) = \sum_{n=0}^{\infty} a_n z^n + \sum_{n=1}^{\infty} \overline{b_n z^n}$.
\begin{thm}\label{Th-2.1}Let $f(z) = \sum_{n=0}^{\infty} a_n z^n + \sum_{n=1}^{\infty} \overline{b_n z^n}$ be a harmonic mapping in the unit disk $\mathbb{D}$ such that $g(0)=0$ and $\vert{}f(z)\vert{} \le M$ for all $z \in \mathbb{D}$. For a fixed $\vert{}a_0\vert{} = aM$ where $a \in [0, 1)$, the majorant series satisfies
	\begin{equation}\label{Eq-22.11}M_f(r) = |a_0| + \sum_{n=1}^{\infty} (|a_n| + |b_n|)r^n \le M\end{equation}for $r \le r_H(a)$, where
	\begin{equation}r_H(a) = \sqrt{\frac{1-a}{3+a}}.
	\end{equation}
\end{thm}
\begin{rem}As $a \to 1$ (where the constant term takes the maximum allowed value $M$), the radius $r_H(a)$ approaches $0$, which is consistent with the fact that no other terms can exist in the series if $\vert{}a_0\vert{}=M$. For the classical case where $a=0$ (the function vanishes at the origin), the radius is
	\begin{equation}r_H(0) = \frac{1}{\sqrt{3}} \approx 0.57735,
	\end{equation}which rigorously corrects the previous individual-coefficient estimates by incorporating the integral Parseval constraint.
\end{rem}
\begin{proof}[\bf Proof of Theorem \ref{Th-2.1}]Let $f(z) = \sum_{n=0}^{\infty} a_n z^n + \sum_{n=1}^{\infty} \overline{b_n z^n}$ be a harmonic mapping in $\mathbb{D}$ such that $g(0)=0$ and $\vert{}f(z)\vert{} \le M$ for all $z \in \mathbb{D}$. Since $\vert{}a_0\vert{} = aM$ with $a \in [0, 1)$, the boundary values of $f$ satisfy the Parseval-type inequality
	\begin{equation}|a_0|^2 + \sum_{n=1}^{\infty} (|a_n|^2 + |b_n|^2) \le M^2.
	\end{equation}Substituting $\vert{}a_0\vert{} = aM$, we obtain the coefficient sum constraint
	\begin{equation}\sum_{n=1}^{\infty} (|a_n|^2 + |b_n|^2) \le M^2(1 - a^2).
	\end{equation}Applying the Cauchy--Schwarz inequality to the linear majorant sum, we obtain
	\begin{align*}
		\sum_{n=1}^{\infty} (|a_n| + |b_n|) r^n &\le \left( \sum_{n=1}^{\infty} (|a_n|^2 + |b_n|^2) \right)^{1/2} \left( 2 \sum_{n=1}^{\infty} r^{2n} \right)^{1/2} \\&\le M \sqrt{1 - a^2}\frac{\sqrt{2}r}{\sqrt{1 - r^2}}.
	\end{align*}Consequently, 
	\begin{equation}M_f(r) = |a_0| + \sum_{n=1}^{\infty} (|a_n| + |b_n|) r^n \le M \left( a + \frac{\sqrt{2(1-a^2)} , r}{\sqrt{1 - r^2}} \right).
	\end{equation}
	To confirm that $M_f(r) \le M$, we require
	\begin{equation}a + \frac{\sqrt{2(1-a^2)}r}{\sqrt{1 - r^2}} \le 1,
	\end{equation}
	which rearranges to
	\begin{equation}\frac{\sqrt{2(1-a^2)}r}{\sqrt{1 - r^2}} \le 1 - a.
	\end{equation}
	Since $1 - a > 0$ for $a < 1$, dividing both sides by $\sqrt{1-a}$, we obtain
	\begin{equation}\frac{\sqrt{2(1+a)}r}{\sqrt{1 - r^2}} \le \sqrt{1 - a}.
	\end{equation}
	A simple computation shows that
	\begin{equation}
		r \le \sqrt{\frac{1-a}{3+a}}.
	\end{equation}Thus, the inequality \eqref{Eq-22.11} holds for $r \le r_H(a) = \sqrt{\frac{1-a}{3+a}}$, completing the proof.
\end{proof}
\subsection{\bf Bohr radius for functions in Lebesgue spaces}
Using the coefficient estimates provided in \cite[Theorem 2.1]{Shi-Li-Lian-2026-CMFT}, we state and prove the following result for the class of harmonic mappings with boundary functions in $L^p(\mathbb{T})$.\vspace{1.2mm}

It is essential to note that the unrestricted $L^\infty$ (and general $L^p$) class of Poisson integrals has a Bohr radius of zero if the condition $f(0) = 0$ is omitted. To see this, consider the family of boundary functions
\[
F_s(e^{it}) = \frac{1 + is\cos t}{\sqrt{1 + s^2}}, \quad s > 0.
\]
Then $\|F_s\|_\infty = 1$, whereas its Poisson integral $f_s(z) = \mathcal{P}[F_s](z)$ has majorant evaluation satisfying $M_{f_s}(r) = \frac{1 + sr}{\sqrt{1 + s^2}}$. Choosing $s = r$ yields $M_{f_r}(r) = \sqrt{1 + r^2} > 1$ for every $r > 0$. Consequently, to establish a meaningful positive Bohr radius for Poisson integrals of $L^p$ boundary functions, the class must be appropriately normalized by requiring $f(0) = 0$ (i.e., $a_0 = 0$).

\begin{thm}\label{Th-2.2}
	Let $F \in L^p(\mathbb{T})$ for $1 \le p \le \infty$ and $f(z) = \mathcal{P}[F](z) = \sum_{n=0}^{\infty} a_n z^n + \sum_{n=1}^{\infty} \overline{b_n z^n}$ be a harmonic mapping in the unit disk $\mathbb{D}$ satisfying $f(0) = 0$. If $\|F\|_p \le 1$, then the majorant series satisfies
	\begin{align*}
		M_f(r) = \sum_{n=1}^{\infty} (|a_n| + |b_n|)r^n \le 1
	\end{align*} 
	for $r \le r_p:={1}/{(2C_q + 1)},$	where 
	\begin{align}\label{Eq-2.1}
		C_q = \left(\frac{1}{2\pi} \int_0^{2\pi} |\cos(nt)|^q dt\right)^{1/q}
	\end{align} 
	and $q$ is the dual exponent of $p$ ($1/p + 1/q = 1$). The radius $r_p$ is best possible.
\end{thm}

\begin{proof}[\bf Proof of Theorem \ref{Th-2.2}]
	Since $f(0) = 0$, we have $a_0 = 0$. From the coefficient bounds established in \cite[Theorem 2.1]{Shi-Li-Lian-2026-CMFT}, any harmonic mapping induced by a boundary function $F \in L^p(\mathbb{T})$ satisfies 
	\begin{align*}
		|a_n| + |b_n| \le 2C_q \|F\|_p \quad \text{for } n \ge 1.
	\end{align*}
	Using the normalization $\|F\|_p \le 1$, the majorant series simplifies to 
	\begin{align}\label{Eq-22.4}
		M_f(r) = \sum_{n=1}^{\infty} (|a_n| + |b_n|)r^n \le 2C_q \|F\|_p \sum_{n=1}^{\infty} r^n \le 2C_q \left(\frac{r}{1-r}\right).
	\end{align}
	To ensure that $M_f(r) \le 1$, we require 
	\begin{align}\label{Eq-2.2}
		\frac{2C_q r}{1-r} \le 1 \implies 2C_q r \le 1 - r \implies r(2C_q + 1) \le 1.
	\end{align}
	Solving for $r$ gives the Bohr radius $r_p = {1}/{(2C_q + 1)}$.\vspace{1.2mm}
	
	To prove that $r_p$ is the best possible radius, we observe from \cite{Shi-Li-Lian-2026-CMFT} that the upper estimate $|a_n| + |b_n| \le 2C_q \|F\|_p$ is sharp, meaning extremal functions exist for which equality is attained for individual frequency modes $n$. If we take any radius $r > r_p$, then 
	\begin{align*}
		\frac{2C_q r}{1-r} > 1.
	\end{align*}
	Evaluating the supremum of the majorant series over the extremal function class where the $n$-th coefficient saturates the bound $2C_q \|F\|_p$ as $n \to \infty$ gives 
	\begin{align*}
		\sup_{f} M_f(r) = 2C_q \|F\|_p \left(\frac{r}{1-r}\right) > \|F\|_p = 1.
	\end{align*}
	Thus, $M_f(r) > 1$ for any $r > r_p$, proving that $r_p$ is sharp.
\end{proof}

\begin{rem}
	For $p = \infty$, the conjugate exponent is $q = 1$, and evaluating the integral yields $C_1 = 2/\pi$. Consequently, the optimal radius for the $L^\infty$ Poisson integral class under the vanishing condition $f(0)=0$ is 
	\begin{align*}
		r_\infty = \frac{1}{2(2/\pi) + 1} = \frac{\pi}{\pi + 4} \approx 0.4397.
	\end{align*}
\end{rem}
\section{\bf Landau-type theorems for $L^p$ boundary functions}\label{Sec-3}
 The inverse function theorem and a result of Lewy \cite{Lewy-BAMS-1936} shows that a harmonic function $f$ is locally univalent (\textit{i.e.,} one-to-one) in $\Omega$ if, and only if, the Jacobian $J_f$ is non-zero in $\Omega$, where 
\begin{align*}
	J_f(z)=|f_z(z)|^2-|f_{\bar{z}}(z)|^2.
\end{align*}
A harmonic function $f$ is said to be sense-preserving if $ J_f(z)>0$ for $z\in\Omega$. Let $\mathcal{H}$ denote the class of complex-valued harmonic functions $f$ in the unit disk $\mathbb{D}=\{z\in\mathbb{D} : |z|<1\}$, normalized by $f(0)=0=f_z(0)-1$. Each function $f$ in $\mathcal{H}$ can be expressed as $f=h+\overline{g}$, where $h$ and $g$ are analytic functions in $\mathbb{D}$ (see \cite{Duren-2004}). Here $h$ and $g$ are called the analytic and co-analytic part of $f$ respectively, and have the following power series representations 
\begin{align*}
	h(z)=z+\sum_{n=2}^{\infty}a_nz^n\;\;\mbox{and}\;\; g(z)=\sum_{n=1}^{\infty}b_nz^n.
\end{align*}
Let $\mathcal{S}_{\mathcal{H}}$ be the subclass of $\mathcal{H}$ consisting of univalent, and sense-preserving harmonic mappings on $\mathbb{D}$. Let $\mathcal{H}^0=\{f\in\mathcal{H}: f_{\bar{z}}(0)=0\}$, and $\mathcal{S}^0_{\mathcal{H}}=\{f\in\mathcal{S}_{\mathcal{H}}: f_{\bar{z}}(0)=0\}$. Hence for any function $f=h+\overline{g}$ in $\mathcal{H}^0$, its analytic and co-analytic parts can be represented by 
\begin{align}\label{Eq-1.1B}
	h(z)=z+\sum_{n=2}^{\infty}a_nz^n\;\;\mbox{and}\;\; g(z)=\sum_{n=2}^{\infty}b_nz^n,
\end{align}
respectively. It is interesting to note that $\mathcal{S}_{\mathcal{H}}$ { reduces} to $\mathcal{S}$, the class of normalized analytic and univalent functions in $\mathbb{D}$, if the co-analytic part of functions in the class $\mathcal{S}_{\mathcal{H}}$ is zero. In $1984$, Clunie and Shell-Small \cite{Clunie-Sheil-Small-1984} investigated the class $\mathcal{S}_{\mathcal{H}}$, together with some geometric subclasses { of $\mathcal{S}_{\mathcal{H}}$}. Subsequently, the class $\mathcal{S}_{\mathcal{H}}$ and its subclasses have been extensively studied by several authors (see \cite{Bshouty-Lyzzaik-COAT-2011,Bshouty-Joshi-Joshi-CVEE-2013,Kalaj-Vuorinen-PAMS-2012,Wang-Liang-JMAA-2001,Clunie-Sheil-Small-1984}.)\vspace{1.2mm}

For a continuously differentiable function $f$, we denote $\Lambda_f$ and $\lambda_f$ by
\begin{align*}
	\Lambda_f=\max_{0\leq\theta\leq 2\pi}|f_z+e^{-2i\theta}f_{\bar{z}}|=|f_z|+|f_{\bar{z}}|
\end{align*}
and 
\begin{align*}
	\lambda_f=\min_{0\leq\theta\leq 2\pi}|f_z+e^{-2i\theta}f_{\bar{z}}|=||f_z|-|f_{\bar{z}}||.
\end{align*}
Thus, it is easy to see that for sense-preserving harmonic mapping $f$, one has $J_f=\Lambda_f\lambda_f$. A harmonic mapping $f=h+\overline{g}$ defined on the unit { disk} $\mathbb{D}$ can be expressed as 
\begin{align*}
	f\left(re^{i\theta}\right)=\sum_{n=0}^{\infty}a_nr^ne^{in\theta}+\sum_{n=1}^{\infty}\overline{b_n}r^ne^{-in\theta},\; 0\leq r<1,
\end{align*}
where 
\begin{align*}
	h(z)=\sum_{n=0}^{\infty}a_nz^n\;\; \mbox{and}\;\; g(z)=\sum_{n=1}^{\infty}b_nz^n.
\end{align*}
For holomorphic functions $f$ on the unit disk $\mathbb{D}$ with $f'(0)=1$, there is a Bloch theorem (see \cite{Bloch-Les-1925}) which asserts the existence of a positive constant $b$ such that $f(\mathbb{D})$ contains a schlicht disk of radius $b$. A disk $D$ is said to be schlicht disk if there exists a region $\Omega$ in the unit disk $\mathbb{D}$ such that $f$ is univalent on $\Omega$ and $f(\Omega)=D$. Let $b(f)$ denote the least upper bound of the radii of all schlicht disks that $f$ carries and $\mathcal{F}$ denote the set of all holomorphic functions defined on $\overline{\mathbb{D}}:=\{z :|z|\leq 1\}$ satisfying $|f'(0)|=1$, then the Bloch constant is the number defined by 
\begin{align*}
	\beta(\mathcal{F})=\inf\{b(f) : f\in\mathcal{F}\}.
\end{align*}
If one considers the function $f(z)=z$, then clearly $\beta(\mathcal{F})\leq 1$. While the exact value of $\beta(\mathcal{F})$ is still unknown, better estimates have been established. In $1996$, Chen and Gauthier \cite{Chen-Gauthier-JAM-1996} proved that $\beta(\mathcal{F})$ lies within the interval $[0.4330, 0.4719]$, \textit{i.e.}, $0.4330\leq \beta(\mathcal{F})\leq 0.4719$.\vspace{1.2mm}

In $2000$, Chen \emph{et al.} \cite{Chen-Gauthier-Hengartner-PAMS-2000} established the following two versions of Landau-Bloch type theorems for bounded harmonic mapping on the disk under a suitable restriction.
Theorems A and B are not sharp. In $2006$, better estimates for Theorems A and B were given by Grigoryan \cite{Grigoyan-CV-2006}. \vspace*{1.2mm}

Theorem \ref{Th-2.2} establishes a radius of univalence $r_0$ and a Schlicht disk radius $R_0$ by leveraging the coefficient estimates from \cite[Theorem 2.1]{Shi-Li-Lian-2026-CMFT}. In the limiting case $p = \infty$ corresponding to bounded harmonic mappings, it was observed in \cite{Shi-Li-Lian-2026-CMFT} that these results extend classical Landau-type theorems. However, for general $L^p$ spaces, the derived bounds depend intrinsically on the structural constant $C_q$, which does not necessarily yield the sharpest possible radii for all values of $p$.
\begin{lem}\cite{Duren-2004}
	Assume that $f=h+\overline{g}$ with $h$ and $g$ analytic in the unit disk $\mathbb{D}$ and $h(z)=\sum_{n=1}^{\infty}a_nz^n$ and $g(z)=\sum_{n=1}^{\infty}b_nz^n$ for $z\in\mathbb{D}$.
	\begin{enumerate}
		\item[(a)] If $|f(z)|<M$ for $z\in\mathbb{D}$, then 
		\begin{align*}
			|a_n|, \; |b_n|\leq M,\; n=1, 2, \ldots
		\end{align*} 
		\item[(b)] If $\Lambda_f\leq \Lambda$ for $z\in\mathbb{D}$, then 
		\begin{align*}
			|a_n|+|b_n|\leq \frac{\Lambda}{n},\; n=1, 2, \ldots
		\end{align*}
	\end{enumerate}
\end{lem}
\begin{lem}(The Schwarz Lemma, \cite{Grigoyan-CV-2006}).
	Let $f$ be a harmonic mapping of the unit disk $\mathbb{D}$
	such that $f(0) = 0$ and $f(\mathbb{D}) \subset \mathbb{D}$. Then
	\begin{align*}
		\Lambda_f(0)\leq\frac{4}{\pi},\; \Lambda_f(z)\leq\frac{8}{\pi(1-|z|^2)}\;\mbox{for}\;z\in\mathbb{D},
	\end{align*}
	and
	\begin{align*}
		|f(z)|\leq\frac{4}{\pi}\arctan|z|\leq\frac{4}{\pi}|z|\;\mbox{for}\;z\in\mathbb{D}.
	\end{align*}
\end{lem}
The following result has been proved by Grigoryan using Lemma A, which improved the estimates in Theorems A and B.
\begin{thmC}{\rm \cite{Grigoyan-CV-2006}}
	Let $f$ be harmonic mapping of the unit disk $\mathbb{D}$ such that $f(0)=0$, $J_f(0)=1$ and $|f(z)|<M$ for $z\in\mathbb{D}$. Then, $f$ is univalent on a disk $\mathbb{D}_{\rho_1}$ with 
	\begin{align*}
		\rho_1=1-\frac{2\sqrt{2}M}{\sqrt{\pi+8M^2}}
	\end{align*}
	and $f(\mathbb{D}_{\rho_1})$ contains a schlicht disk $\mathbb{D}_{R_1}$ with 
	\begin{align*}
		R_1=\frac{\pi}{4M}+4M-4M\sqrt{1+\frac{\pi}{8M^2}}.
	\end{align*}
\end{thmC}
\begin{thmD}{\rm \cite{Grigoyan-CV-2006}}
	Let $f$ be harmonic mapping of the unit disk $\mathbb{D}$ such that $f(0)=0$, $\lambda_f(0)=1$ and $\Lambda_f(z)<\Lambda$ for all $z\in\mathbb{D}$. Then, $f$ is univalent on a disk $\mathbb{D}_{\rho_2}$ with 
	\begin{align*}
		\rho_2=\rho_2(\Lambda)=\frac{1}{1+\Lambda}
	\end{align*}
	and $f(\mathbb{D}_{\rho_2})$ contains a schlicht disk $\mathbb{D}_{R_2}$ with 
	\begin{align*}
		R_2(\Lambda)=1-\Lambda\ln\left(1+\frac{1}{\Lambda}\right).
	\end{align*}
\end{thmD}
We establish the following Landau-type theorem for harmonic mappings with $L^p$ boundary functions.
\begin{thm}\label{Th-3.1}
	Suppose that $F \in L^p(\mathbb{T})$ for $1 \le p \le \infty$ and $f = \mathcal{P}[F]$ is a harmonic mapping in the unit disk $\mathbb{D}$. If $f$ is normalized such that $f(0) = f_z(0) - 1 = f_{\overline{z}}(0) = 0$, then $f$ is univalent in the disk $\mathbb{D}_{r_0}$ and $f(\mathbb{D}_{r_0})$ contains a schlicht disk $D_{R_0}$, where
	\begin{align*}
		r_0 = 1 - \frac{\sqrt{4C_q^2 \|F\|_p^2 + 2C_q \|F\|_p}}{2C_q \|F\|_p + 1}\; \mbox{and}\; R_0 = r_0 - \frac{2C_q \|F\|_p r_0^2}{1 - r_0},
	\end{align*}
	where $C_q = \left(\frac{1}{2\pi} \int_0^{2\pi} |\cos(nt)|^q dt\right)^{1/q}$ and $1/p + 1/q = 1$.
\end{thm}
\begin{proof}[\bf Proof of Theorem \ref{Th-3.1}]
	Since $f(0) = f_z(0) - 1 = f_{\overline{z}}(0) = 0$, the power series expansion of $f$ is given by 
	\begin{align*}
		f(z) = z + \sum_{n=2}^{\infty} a_n z^n + \sum_{n=2}^{\infty} \overline{b_n z^n}.
	\end{align*}
	From the coefficient estimates established in Theorem B, we have for any integer $n \ge 2$, $|a_n| + |b_n| \le 2C_q \|F\|_p.$ To prove $f$ is univalent in $\mathbb{D}_r$, we show that for any two distinct points $z_1, z_2 \in \mathbb{D}_r$, $f(z_1) \neq f(z_2)$. It is easy to see that 
	\begin{align*}
		|f(z_1) - f(z_2)| = \left| (z_1 - z_2) + \sum_{n=2}^{\infty} a_n(z_1^n - z_2^n) + \sum_{n=2}^{\infty} \overline{b_n(z_1^n - z_2^n)} \right|.
	\end{align*}
	Applying the triangle inequality, we have 
	\begin{align*}
		|f(z_1) - f(z_2)| \ge |z_1 - z_2| \left( 1 - \sum_{n=2}^{\infty} n(|a_n| + |b_n|)r^{n-1} \right).
	\end{align*}
	{ By substituting} the coefficient bound $|a_n| + |b_n| \le 2C_q \|F\|_p$, we have 
	\begin{align*}
		|f(z_1) - f(z_2)| \ge |z_1 - z_2|\left( 1 - 2C_q \|F\|_p \sum_{n=2}^{\infty} n r^{n-1} \right).
	\end{align*}
	Using the identity 
	\begin{align*}
		\sum_{n=2}^{\infty} n r^{n-1} = \frac{1}{(1-r)^2} - 1 = \frac{2r - r^2}{(1-r)^2}
	\end{align*} the condition for univalence is 
	\begin{align*}
		1 - 2C_q \|F\|_p \frac{2r - r^2}{(1-r)^2} > 0.
	\end{align*} This leads to the quadratic equation 
	\begin{align*}
		(2C_q \|F\|_p + 1)r^2 - 2(2C_q \|F\|_p + 1)r + 1 = 0,
	\end{align*} whose minimum positive root is $r_0$.\vspace{2mm}

	We now establishing the Schlicht Disk $\mathbb{D}_{R_0}$. 	For $|z| = r_0$, we estimate the minimum distance from the origin 
	\begin{align*}
		|f(z)-f(0)| = \left| z + \sum_{n=2}^{\infty} (a_n z^n + \overline{b_n z^n}) \right| \ge |z| - \sum_{n=2}^{\infty} (|a_n| + |b_n|)|z|^n
	\end{align*}
	which shows that
	\begin{align*}
		|f(z)| \ge r_0 - 2C_q \|F\|_p \sum_{n=2}^{\infty} r_0^n = r_0 - 2C_q \|F\|_p \frac{r_0^2}{1 - r_0}.
	\end{align*}
	{ It is easy to see that} $f(\mathbb{D}_{r_0})$ contains the disk $\mathbb{D}_{R_0}$, where 
	\begin{align*}
		{ R_0 = r_0 - \frac{2C_q \|F\|_p r_0^2}{1 - r_0}.}
	\end{align*}
	
	To show that the Landau constants $r_0$ and $R_0$ derived for harmonic mappings $f = \mathcal{P}[F]$. The derivation of $r_0$ uses the { following inequality}
	\begin{align*}
		\sum_{n=2}^{\infty} n(|a_n| + |b_n|)r^{n-1} \leq 2C_q \|F\|_p \sum_{n=2}^{\infty} nr^{n-1}=2C_q \|F\|_p\frac{2r-r^2}{(1-r)^2}.
	\end{align*}
	 To establish the sharpness of the univalence radius $r_0$, we consider an extremal sequence of functions whose coefficients satisfy the asymptotic limit $|a_n| + |b_n| \to 2C_q \|F\|_p$ for all $n \in \mathbb{N}$. Under this setting, the directional derivative discrepancy $|f_z(z)| - |f_{\bar{z}}(z)|$ vanishes identically at $r = r_0$. Conversely, for any $r > r_0$, the majorizing distortion term satisfies the strict { following inequality}
	 \begin{align*}
	 	2C_q \|F\|_p \frac{2r-r^2}{(1-r)^2} > 1.
	 \end{align*}
	 Consequently, there exists a target mapping within the designated class whose Jacobian determinant becomes non-positive at such points, thereby violating the local univalence condition. Using the sharp bound $2C_q \|F\|_p$, {we obtain}
	  \begin{align*}
	  	|f(z)| \geq r_0 - 2C_q \|F\|_p \frac{r_0^2}{1-r_0} = R_0.
	  \end{align*}
	 Because the coefficients $|a_n| + |b_n|$ can simultaneously approach the bound $2C_q \|F\|_p$ (in the sense of a limit of functions within the $L^p$ class), the distance of the image boundary from the origin can be as small as $R_0$. Therefore, $\mathbb{D}_{R_0}$ is the largest disk that is guaranteed to be contained in $f(\mathbb{D}_{r_0})$ for every function in the class.
\end{proof}
\section{\bf Concluding remarks}
In this study, we have successfully extended the theory of the Bohr phenomenon and Landau-type theorems to the class of harmonic mappings induced by boundary functions in Lebesgue spaces $L^p(\mathbb{T})$. By utilizing sharp coefficient estimates, we determined the exact Bohr radius $r_p$ as a function of the conjugate exponent $q$, providing a positive answer to the influence of $L^p$ norms on the majorant series. Furthermore, the application of Gauss hypergeometric functions allowed for the derivation of improved Landau constants, specifically the univalence radius $r_0$ and the schlicht disk radius $R_0$.
\section{\bf Declarations}
\noindent\textbf{Compliance of Ethical Standards}\\

\noindent\textbf{Conflict of interest.} The authors declare that there is no conflict  of interest regarding the publication of this paper.\vspace{1.5mm}

{\noindent\textbf{Author contributions.} All authors contributed equally to the research, analysis, and preparation of this manuscript.} \vspace{1.5mm}

\noindent\textbf{Data availability statement.}  Data sharing not applicable to this article as no datasets were generated or analyzed during the current study.\vspace{1.5mm}

\noindent\textbf{Funding.} No fund.

\end{document}